\documentclass[review]{elsarticle}

\usepackage[linktocpage=true,pagebackref=true,linkcolor=black,citecolor=black]{hyperref}

\journal{Journal of \LaTeX\ Templates}

\usepackage{verbatim}
\usepackage[utf8]{inputenc}
\usepackage[english]{babel}
\usepackage{amssymb}
\usepackage{amsthm}
\usepackage{amsmath,amssymb,amsfonts,graphicx,float,multicol,fleqn}
\usepackage{subcaption}
\usepackage{multirow}

\UseRawInputEncoding

\usepackage{algorithm}
\usepackage[noend]{algpseudocode}

\newcommand{\IR}{\ensuremath{\mathbb{R}}} 

\usepackage[letterpaper, margin=1in]{geometry}
\usepackage[dvipsnames]{xcolor}
\usepackage{enumerate}

\usepackage{caption} 
\theoremstyle{definition}
\newtheorem{example}{Example}[section]

\newcommand{\IN}{\ensuremath{\mathbb{N}}}

\begin{document}

\begin{frontmatter}

\title{JDNN: Jacobi Deep Neural Network for Solving Telegraph Equation}

\author[1]{Maryam Babaei}
\ead{mar_babaei@sbu.ac.ir}
\author[1]{Kimia Mohammadi Mohammadi}
\ead{kimia.moh97@gmail.com}
\author[1]{Zeinab Hajimohammadi}
\ead{Z\_Hajimohammadi@sbu.ac.ir}
\author[1,3]{Kourosh  Parand
\corref{mycorrespondingauthor}}
\cortext[mycorrespondingauthor]{Corresponding author}
\address[1]{Department of Computer and Data Sciences, Faculty of Mathematical Sciences, Shahid Beheshti University, Tehran, Iran}
\address[3]{Institute for Cognitive and Brain Sciences, Shahid Beheshti University, Tehran, Iran}

\begin{abstract}
In this article, a new deep learning architecture, named JDNN, has been proposed to approximate a numerical solution to Partial Differential Equations (PDEs). The JDNN is capable of solving high-dimensional equations. Here, Jacobi Deep Neural Network (JDNN) has demonstrated various types of telegraph equations. This model utilizes the orthogonal Jacobi polynomials as the activation function to increase the accuracy and stability of the method for solving partial differential equations. The finite difference time discretization technique is used to overcome the computational complexity of the given equation. The proposed scheme utilizes a Graphics Processing Unit (GPU) to accelerate the learning process by taking advantage of the neural network platforms. Comparing the existing methods, the numerical experiments show that the proposed approach can efficiently learn the dynamics of the physical problem. \end{abstract}


\end{frontmatter}

\section{Introduction}
Using partial differential equations (PDEs) in modeling various types of physical phenomena is prevalent. Oliver Heaviside introduced the telegraph equation, one of the most critical problems studied over the last decades. Generally, it has been used to model the vibrations of structures and is the basis for the fundamental equations of atomic physics \cite{par}.

In the present work, we are dealing with the following second-order hyperbolic problem:

\begin{equation}\label{qlm1}
\hspace*{\fill} \frac{\partial ^{2}u}{\partial t^{2}}+ 2\alpha \frac{\partial u}{\partial t} + \beta^2 u=\nabla^2u+f(\overline{x},t),\hspace*{\fill}
\end{equation}
where $\overline{x}$ is a $d$-dimensional vector, $u$ is a function of $\overline{x}$ and time variable $t$ and $\nabla$ is the gradient operator. Eq.  (\ref{qlm1}), known as the second-order $d$-dimensional telegraph equation with constant coefficients, is commonly used in signal analysis for transmission and propagation of electrical signals, along with its applications in other fields \cite{ind}.
Researchers have studied numerical schemes such as finite differences, spectral methods, and finite elements to achieve solutions for the telegraph equation \cite{ar44, articl3e, irk2019numerical}. Some of these mesh-based methods are applied to spatial discretization, such as the Finite Volume Method (FVM) and Boundary Element Method (BEM) \cite{doi:10.2514/3.559, Wrobel2003BoundaryEM}. Another group of numerical methods is known as meshless methods that do not require mesh for discretization, like Element Free Galerkin (EFG) \cite{https://doi.org/10} and Meshless Local Petrov-Galerkin (MLPG) \cite{article3, article1}. Jiwari has proposed two differential quadrature methods based on the Lagrange interpolation and modified cubic B-splines to find the approximate solution of one and two-dimensional hyperbolic equations, such as the telegraph equation \cite{artJ}. Oru{\c{c}} has introduced an algorithm based on the Hermite wavelets for the 2-dimensional hyperbolic telegraph equation \cite{arO}. Jiwari et al. have presented a numerical technique based on the polynomial differential quadrature method for the second-order 1-dimensional hyperbolic telegraph equation \cite{artiji}.
Recently, the integration of machine learning and deep learning with numerical methods has been proposed to solve PDE equations \cite{hajimohammadi2021legendre, 6224185}.

Machine Learning (ML) consists of computer algorithms that analyze data to improve their accuracy. As a subset of machine learning, deep learning aims to model human brain behavior. It provides reason-able solutions for various types of applications, such as image recognition and natural language processing \cite{GoodBengCour16}.
The machine learning approaches have been leveraged in mathematical applications based on previous records. Mehrkanoon and Suykens have proposed an approach based on Least Squares Support Vector Regression (LS-SVR) for solving second-order PDEs with variable coefficients \cite{MEHR}. 
They have shown that using these methods can result in an outstanding performance. Consequently, researchers have developed deep neural networks to learn the physical dynamics of a wide range of engineering problems \cite{SIRI, RAISI}.

On the other side of science, orthogonal functions have been utilized in approximation theory and numerical analysis.
These functions' applications also appeared
in kernel-based machine learning methods. To name a few, Padierna et al. have developed Gegenbauer polynomials in Support Vector Machines (SVMs) to classify some real-world problems \cite{padierna2018novel}. Wavelet kernels have been introduced by Zhang et al. to improve the accuracy of SVMs \cite{wavelet2}. Gupta et al. have evaluated the applications of these functions in principal component analysis (PCA) \cite{wavelet}.
The orthogonal polynomials have also been used in LS-SVMs to solve integrals and differential equations \cite{HAJI, parand2021parallel, parand2}.
Additionally, they can be used in the structure of neural networks to improve their performance.  For instance, a single layer fractional orthogonal neural network with fractional order of Legendre functions has been presented by Hadian et al. \cite{HADI}.
Hajimohamadi et al. have recently extended this idea to the deep neural network.\cite{hajim}
Chakraverty and Mall introduced artificial neural networks with orthogonal polynomials as activation functions \cite{chak}.

In the current paper, a new neural network architecture, which will be called JDNN, 
is designed to increase the model's accuracy in simulating telegraph equations.  In order to achieve this goal, finite difference and spectral methods are utilized to improve efficiency.

The rest of the paper is organized as follows: In section 2, we will explain different types of telegraph equations along with their initial and boundary conditions. We will discuss Jacobi, Legendre, and Chebyshev polynomials in the third section. In addition, we will briefly discuss exciting finite difference methods for the telegraph equation. In section 4, we will describe the JDNN method. This approach helps us to solve differential equations in both continuous and discrete-time forms. The construction and verification of JDNN will be illustrated in Section 5. At last, we will give the concluding remarks in Section 6.
\section{Mathematical formulations of the telegraph equation}
The telegraph equation is a PDE That has been used in many fields like transmitting digital and analog signals \cite{doi:10.1063/1.369258}, the random walk theory \cite{article22}, microwaves, and radiofrequency fields \cite{Roussy1998FoundationsAI}. Some numerical and analytical solution approaches were proposed for this problem. The history of research can be seen in \cite{article44, MA2016236, doi:10.1080/00207160211918}.

In this article, we consider the one and two-dimensional telegraph equation. In a 1-dimensional case, it can be expressed as:

\begin{equation}\label{qlm1m}
\hspace*{\fill} \frac{\partial ^{2}u}{\partial t^{2}}+ 2\alpha \frac{\partial u}{\partial t} + \beta u=\frac{\partial ^{2}u}{\partial x^{2}}+f(x,t), \quad x \in \Omega \subset \mathbb{R}  ,t>0, \hspace*{\fill}
\end{equation}
where $\alpha$ and $\beta$ are constant coefficients and $\Omega \in [A,B]\subset \mathbb{R}$. The initial and boundary conditions are considered in the following form:
\begin{equation}\label{qlm2m}
\hspace*{\fill} u(x,0)=g_1(x),\quad x \in \Omega , \hspace*{\fill}
\end{equation}
\begin{equation}\label{qlm3m}
\hspace*{\fill} \frac{\partial u}{\partial t}(x,0)=g_2(x),\quad x \in \Omega ,  \hspace*{\fill}
\end{equation}
\begin{equation}\label{qlm4m}
\hspace*{\fill} u(x,t)=h_1(x),\quad x =A ,t>0, \hspace*{\fill}
\end{equation}
\begin{equation}\label{qlmmm}
\hspace*{\fill} u(x,t)=h_2(x),\quad x =B ,t>0. \hspace*{\fill}
\end{equation}
The two-dimensional telegraph equation is defined as
\begin{equation}\label{qlm1n}
\hspace*{\fill} \frac{\partial ^{2}u}{\partial t^{2}}+ 2\alpha \frac{\partial u}{\partial t} + \beta^2 u=\nabla^2u+f(x_1,x_2,t), \quad x_1,x_2 \in \Omega \subset \mathbb{R}^2  ,t>0, \hspace*{\fill}
\end{equation}
where $\nabla=<\frac{\partial}{\partial x_1}, \frac{\partial}{\partial x_2}>$. Initial and boundary conditions are as follows: 
\begin{equation}\label{qlm2n}
\hspace*{\fill} u(x_1,x_2,0)=g_1(x_1,x_2),\quad \frac{\partial u}{\partial t}(x_1,x_2,0)=g_2(x_1,x_2),\quad x_1, x_2 \in \Omega , \hspace*{\fill}
\end{equation}
\begin{equation}\label{qlm3n}
\hspace*{\fill} u(x_1,x_2,t)=h_1(x_1,x_2),\quad u(x_1,x_2,t)=h_2(x_1,x_2),\quad x_1, x_2\in \Gamma. \hspace*{\fill}
\end{equation}
Here, $\Gamma$ is the boundary points of $\Omega$.
\section{Preliminaries}

We will first discuss the Jacobi polynomials and the finite difference method before moving to the main architecture.


\subsection{Jacobi polynomials}
The Jacobi polynomials are a family of the orthogonal polynomials denoted by $J_n^{\alpha, \beta}(x)$ with parameters $\alpha,\beta \in \IR ^{>-1}$ of the nth degree \cite{funaro2008polynomial}. Jacobi orthogonal polynomials are defined on the interval $[-1,1]$ with respect to the weight function $\omega^{\alpha, \beta}(x)=(1-x)^{\alpha}(1+x)^{\beta}$. They can be defined in both explicit and implicit expressions.
The singular Sturm-Liouville problem for Jacobi polynomials is characterized by \cite{funaro2008polynomial}:
\begin{equation}
    -(1-x^2)\frac{d^2}{\partial x^2}J_n^{\alpha, \beta}(x)+((\alpha+\beta+2)x+\alpha-\beta)\frac{d}{\partial x}J_n^{\alpha, \beta}(x)-(n(n+\alpha+\beta+1))J_n^{\alpha, \beta}(x)=0\quad n \in \IN.
\end{equation}
The orthogonality of their orthogonality is determined in the following form \cite{funaro2008polynomial}:\\
\begin{equation}
    \hspace*{\fill}\int_{-1}^{1} J_n^{\alpha, \beta}(x)J_m^{\alpha, \beta}(x)\omega^{\alpha, \beta}(x),dx =\gamma_n^{\alpha, \beta}\delta_{n,m}.\hspace*{\fill}
\end{equation}
where $\gamma_n^{\alpha, \beta}=	\|	J_n^{\alpha, \beta}\|_{\omega^{\alpha, \beta}}^2$ and $\delta_{n,m}$ is  the Kronecker delta function.\\
The recursive formula for $\{ J_n^{\alpha, \beta}(x)\}_{n \geq 0}$ becomes as follows \cite{funaro2008polynomial}:\\
\begin{eqnarray}
&& J_0^{(\alpha,\beta)}(x) = 1,\nonumber\\
&& J_1^{(\alpha,\beta)}(x) =\frac{1}{2}(\alpha + \beta +2)x + \frac{1}{2}(\alpha - \beta),\nonumber\\
&& J_n^{(\alpha,\beta)}(x) =(\rho_n x+\sigma_n)J_{n-1}^{(\alpha,\beta)}(x)+\tau_n J_{n-2}^{(\alpha,\beta)}(x), \quad n \geq 2,\\
&& \rho_n = \frac{(2n + \alpha + \beta)(2n + \alpha + \beta -1)}{2n(n+ \alpha + \beta)},\nonumber\\
&&\sigma _n = \frac{(\alpha ^2 - \beta ^2)(2n + \alpha + \beta -1)}{2n(n+ \alpha + \beta)(2n + \alpha + \beta -2)},\nonumber\\
&& \tau_n = -\frac{(n+ \alpha -1)(n + \beta -1 )(2n + \alpha + \beta)}{n(n+\alpha + \beta)(2n + \alpha +\beta -2)}\nonumber.
\end{eqnarray}
Choosing a proper set of Jacobi polynomials parameters leads to obtaining the other well-known polynomials, including Chebyshev and Legendre polynomials. Legendre polynomials are obtained with $\alpha=\beta=0$, which is denoted by $L_n(x)$, and have the following recurrence relation \cite{funaro2008polynomial}:\\
\begin{eqnarray}
&& L_0 (x)=1, \quad L_1 (x)=x, \nonumber\\
&& (n+1) L_{n+1}(x)=(2n+1)x L_n(x)-n L_{n-1}(x),\quad n\ge 1.
\end{eqnarray}
Also, Legendre polynomials have orthogonality with the weight function $\omega(x)=1$ in the $[-1,1]$ domain. Furthermore, the first kind of Chebyshev polynomials is obtained with $\alpha=\beta=-\frac{1}{2}$, which is denoted by $T_n(x)$, and has the following recurrence relation \cite{funaro2008polynomial}:\\
\begin{eqnarray}
&& T_0 (x)=1, \quad T_1 (x)=x,\nonumber\\
&& T_{n+1}(x)=2x T_n(x)- T_{n-1}(x),\quad n\ge 1.
\end{eqnarray}
The weight function for the first kind of Chebyshev polynomials is $\omega(x)=\frac{1}{\sqrt{1-x^2}}$  in the corresponding domain.
\subsection{Finite difference methods}
Finite difference methods are a class of numerical procedures for approximating the solution of mathematical equations. Taking advantage of the discretization techniques, they simplify the approximation process and reduce its computational complexity. 
In this section, we will first recall the finite difference approximation of the derivative. Then, we will apply numerical differentiation to the telegraph equation.

A first-order finite difference approximation $U'(t)=\frac{\partial u}{\partial t}$ With the error denoted by $o(h)$ error can be estimated by Taylor\textquotesingle s theorem \cite{johnson2012numerical}:
\begin{equation}\label{qlm6}
\hspace*{\fill} U'(t)\approx \frac{U(t+h)-U(t)}{h},\hspace*{\fill}
\end{equation}
where $h$ is a step size, additionally, the second finite difference approximation $U''(t)$  with $O(h)$ error is defined as:
\begin{equation}\label{qlm10}
\hspace*{\fill} U''(t) \approx \frac{U(t-h)-2U(t)+U(t+h)}{h^2},\hspace*{\fill}
\end{equation}
Now we can write Eq. (\ref{qlm1}) as follow:
\begin{equation}\label{qlm11}
\hspace*{\fill} \frac{U_{i-1}-2U_{i}+U_{i+1}}{h^2}+ 2\alpha \frac{U_{i+1}-U_{i}}{h}= -\beta^2 U_{i+1}+\frac{\partial ^{2}U_{i+1}}{\partial x^{2}}+f_{i+1},\hspace*{\fill}
\end{equation}
Where $U_{i+n}=U(\overline{x},t+nh)$, $f_{i+n}=f(\overline{x},t+nh)$, and $\overline{x}$ is a d-dimension vector. It is an approximation method for the telegraph equation used in our network that will be discussed in the next section.

In this part, we consider $h=\Delta t$ and the left hand of Eq. (\ref{qlm11}) as:
\begin{equation}
\hspace*{\fill}N(u(\overline{x},t+\Delta t))= -\beta^2 U(\overline{x},t+h)+\frac{\partial^{2}U(\overline{x},t+h)}{\partial x^{2}}+f(\overline{x},t+h).\hspace*{\fill}
\end{equation}
Then we can define the residual function for our method as follows:
\begin{equation}\label{qlm1236}
\hspace*{\fill}Res(\overline{x},t+\Delta t)=(1+2\alpha\Delta t)U(\overline{x},t+\Delta t)-2(1+\alpha \Delta t)U(\overline{x},t)+U(\overline{x},t-\Delta t)-\Delta t^2N(u(\overline{x},t+\Delta t)).\hspace*{\fill}
\end{equation}

We set $\hat{u}_{t_3}=U(\overline{x},t+\Delta t)$, $u_{t_2}=U(\overline{x},t)$, and $u_{t_1}=U(\overline{x},t-\Delta t)$ so that $\hat{u}_{t_3}$ is the output of the JDNN and the final output at $t+\Delta t$ is illustrated by $u^*_{t_3}$. The goal of this approach is to minimize residual function.
\section{Description of JDNN}
In this section, we will demonstrate the Jacobi Deep Neural Network (JDNN) construction to approximate the telegraph equation's solution.

JDNN architecture contains two networks. The first part is a deep multilayer neural network with the initial layer of Jacobi polynomials as activation functions. This network with n layers is defined as follows:
\begin{eqnarray}
&& \mathcal{H}_0=2(\overline{x} - min(\overline{x})) / (max(\overline{x}) - min(\overline{x})) -1 , \quad \overline{x} \in \mathbb{R}^d,\nonumber\\
&& \mathcal{H}_1=J(W^{(1)}\mathcal{H}_0+b^{(1)}),\nonumber\\
&&\mathcal{H}_i=Tanh(W^{(i)}\mathcal{H}_{i-1}+b^{(i)}) \quad 2\leq i\leq n-1,\\
&&\mathcal{H}_n=W^{(n)}\mathcal{H}_{n-1}+b^{(n)}\nonumber.
\end{eqnarray}
Where $W^{(i)}$ and $b^{(i)}$ for $i=0,1,...,n$ are the weight and the bias parameters. $\mathcal{H}_0$ is the input layer with a d-dimension normalized in the range $[-1, 1]$.
$J=[J_1,...,J_n]^{T}$ in the first hidden layer includes an array of different degrees of Jacobi orthogonal polynomials, and $\mathcal{H}_1$ is named an orthogonal layer. Choosing appropriate $\alpha$ and $\beta$ in the Jacobi polynomial formula can use Legendre and Chebyshev polynomials as activation functions.

Moreover, the hyperbolic tangent function (Tanh) is used in other hidden layers. Using Tanh approximates the answer space to a more complex space than polynomials. $\mathcal{H}_n$ is the output layer that gives the value of the function.

\begin{figure}[H]
	\centering
	\includegraphics[width=0.95\linewidth]{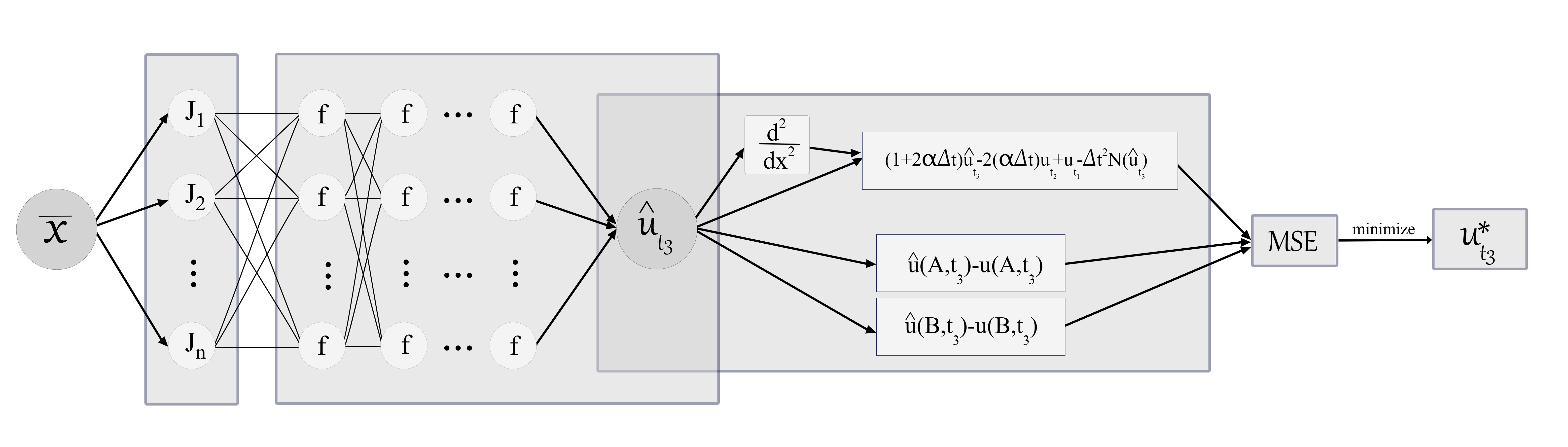}
	\caption{Schematic of a Jacobi deep neural network (JDNN), $J_n$ nodes are the different degrees of the Jacobi polynomials as activation functions, and f nodes contain Tanh activation function. $\hat{u}_{t3}$ is the output of the first part in each epoch, and $u^*_{t3}$ is the overall best result. The second part of JDNN Includes operating nodes for calculating the residual form of the telegraph equation and boundary conditions of $u$ in the Target time step.}
	\label{fig}
\end{figure}

In the second network, the residual form of the telegraph equation is obtained by applying the operating nodes according to the residual function of the telegraph equation. We have used automatic differentiation for derivative calculations that we need for calculating the residual form. Applying automatic differentiation allows derivatives to be performed quickly and accurately. The network loss function is constructed according to the residual form of the equation and boundary conditions. The loss function of the network is as follows:
\begin{equation}\label{loss}
\hspace*{\fill}Loss=MSE_{Res}+MSE_{BC},\hspace*{\fill}
\end{equation}
where MSE refers to the mean squared error, which is calculated in the following form:
\begin{equation}
\hspace*{\fill}MSE_{Res}=\frac{1}{N_{Res}} \sum_{i=1}^{N_{Res}} {|\hat{f_{Res}}(\overline{x}^{i}),t_{3}|}^2 ,\hspace*{\fill}
\label{mse1}
\end{equation}
\begin{equation}
\hspace*{\fill}MSE_{BC}=\frac{1}{N_{BC}} \sum_{i=1}^{N_{BC}} {|\hat{u}(A, t_{3})-u(A, t_{3})|}^2+{|\hat{u}(B, t_{3})-u(B, t_{3})|}^2 .\hspace*{\fill}
\label{mse2}
\end{equation}
Where A and B are boundary points of the equation, $ N_{Res}$ and $N_{BC}$ are the numbers of training points in the range $[A, B]$ and training points for boundary conditions, respectively. Loss minimization is obtained by performing the Adam algorithm\cite{adam}, and a Quasi-Newton Method called L-BFGS-B \cite{bfgs}. The maximum number of iterations is specified to stop the Adam algorithm, But the L-BFGS-B method is stopped when it converges. Figure \ref{fig} shows the architecture of JDNN for solving the telegraph equation.
\section{Numerical results}
In this paper, we have solved four examples of the telegraph equation by applying JDNN. To present the JDNN results, we considered a similar neural network without the orthogonal function (simple DNN). Numerical results for each of them Have been investigated in this section.

We use four types of errors: $L_2$, relative $L_2$, $L_{\infty}$, and RMS error which are calculated in the following form Respectively:
\begin{eqnarray}
&& L_2 \; error= \|u-\hat{u}\|_2=\sqrt{\sum_{i}|u_i-\hat{u}_i|^2},\nonumber\\
&&  relative\; L_2 \; error= \frac{\|u-\hat{u}\|_2}{\|u\|_2},\nonumber\\
&& L_{\infty} \; error= Max_i|u_i-\hat{u}_i|,\nonumber\\
&& RMS \; error= \sqrt{\sum_{i}\frac{|u_i-\hat{u}_i|^2}{n}}\nonumber.
\end{eqnarray}
To show the JDNN architecture in the tables in a better way, we display the layers of each network as $[I, J, F_1, F_2, ..., F_n]$. Where $I$ is the diminution of $\overline{x}$ and $J$ is the number of neurons with orthogonal polynomials activation functions, $F_i$ for $2\leq i\leq n-1$ is the number of neurons with Tanh activation function, and $F_n$ is the output layer. In addition, we Have used the maximum iteration of the Adam algorithm up to 5000 times.
\begin{example}\label{ex1}
Considering
$\alpha=\frac{1}{2}$
and
$\beta=1$
at the interval 
$0\leq x\leq 4$
in hyperbolic telegraph Eq. (\ref{qlm1}) with the following initial conditions:
\begin{equation*}
\hspace*{\fill}
\left\{
\begin{array}{ll}
u(x,0)=e^x & 0\leq x\leq 4\\
u_t(x,0)=-e^x, & 0\leq x\leq4
\end{array}
\right.
\hspace*{\fill}
\end{equation*}
we have
\begin{equation*}
\hspace*{\fill} \frac{\partial ^{2}u}{\partial t^{2}}+ \frac{\partial u}{\partial t} + u=\frac{\partial ^{2}u}{\partial x^{2}},\hspace*{\fill}
\end{equation*}
The analytical solution is shown in \cite{MOMANI20051126} as follows:
\begin{equation*}
\hspace*{\fill}u(x,t)=exp(x-t).\hspace*{\fill}
\end{equation*}
In this example, we use seven hidden layers with 8 Legendre polynomials as activation functions in the first. Evaluating function values in time $t_1=0.6$ and $t_2=0.8$, our proposed model approximates function for $t_3=1$. It is trained and tested with the number of random points equal to 200 and 100, respectively. Moreover, in figure \ref{fig2}, the numerical solutions and the residual
function on the interval $[0, 4]$ for $t=1$ have been presented. Figure \ref{fig3} compares the trend of loss function changes for these two networks, showing the proposed method's stability.
In table  \ref{1.1}, the relative $L_2$ errors for JDNN and the simple DNN have been reported. Besides, Table \ref{2.1}, compares the results with the collocation points method presented by Dehghan et al. \cite{https://doi.org/10.1002/num.20306}.
\begin{figure}[H]
	\centering
	\begin{subfigure}[b]{0.4\textwidth}
      \includegraphics[width=1\linewidth]{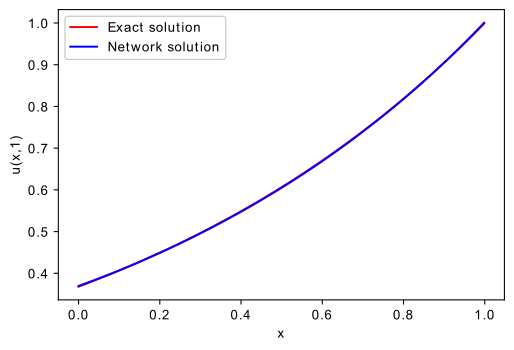}
      \caption{}
    \end{subfigure}
    \begin{subfigure}[b]{0.4\textwidth}
      \includegraphics[width=1.1\linewidth]{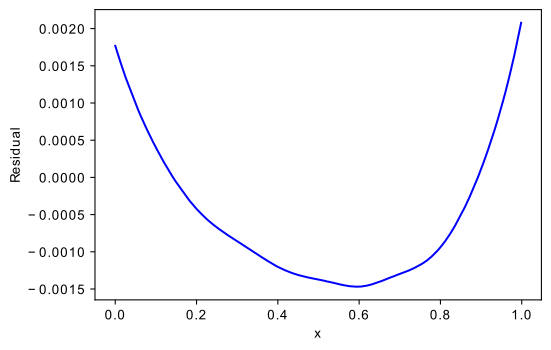}
      \caption{}
    \end{subfigure}
	\caption{(a) Exact and JDNN solutions for the example \ref{ex1} with training points in [0, 4] and (b) Absolute of the residual function.}
        \label{fig2}
\end{figure}
\begin{figure}[H]

	\centering
	\begin{subfigure}[b]{0.4\textwidth}
  	\includegraphics[width=1.0\linewidth]{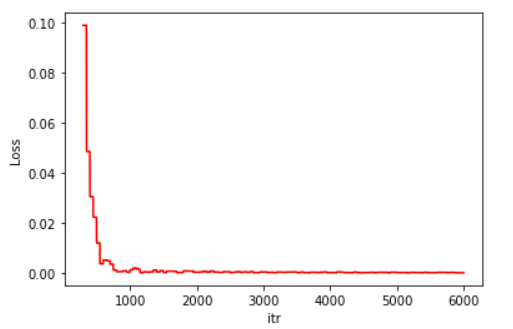}
      \caption{}
    \end{subfigure}
    \begin{subfigure}[b]{0.4\textwidth}
  	\includegraphics[width=1.0\linewidth]{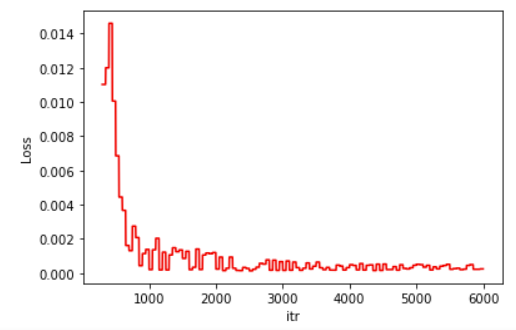}
      \caption{}
    \end{subfigure}
	\caption{\label{fig3}Comparing the trend of Loss function changes in (a) the JDNN and (b) the simple DNN.}

\end{figure}
\begin{table*}[!hbtp]
	\centering
	\caption{Comparison of relative $L_2$ errors obtained in the JDNN and the simple DNN.}
	\begin{tabular}{|c|c|c|c|}
		\hline  \multirow{1}{*}{Method}& \multicolumn{1}{|c|}{Layers}& \multicolumn{1}{|c|}{Train data error}&  \multicolumn{1}{|c|}{Test data error}
		\\
		\hline JDNN & $[1, 8, 16, 16, 32, 16, 16, 8, 1]$ &  $7\times 10^{-4}$ &  $9\times 10^{-4}$\\
		Simple DNN & $[1, 8, 16, 16, 32, 16, 16, 8, 1]$ & $1.5\times 10^{-3}$ &  $1.6\times 10^{-3}$\\
		\hline
	\end{tabular}
	\label{1.1}
\end{table*}

\begin{table*}[!hbtp]
	\centering
	\caption{Comparison of the results obtained in the JDNN with the Collocation points method.}
	\begin{tabular}{|c|c|c|c|}
		\hline  \multirow{1}{*}{Method} &  \multicolumn{1}{|c|}{$L_{\infty}$-error} &  \multicolumn{1}{|c|}{$L_{2}$-error} &  \multicolumn{1}{|c|}{RMS}
		\\
		\hline JDNN  &  $1.8625\times 10^{-5}$&  $8.2376\times 10^{-5}$&  $2.2476\times 10^{-6}$\\
		Collocation points \cite{https://doi.org/10.1002/num.20306}  &  $2.2931\times 10^{-5}$&  $1.7163\times 10^{-4}$&  $8.5711\times 10^{-6}$\\
		\hline
	\end{tabular}
	\label{2.1}
\end{table*}
\end{example}
\begin{example}\label{ex2}
Considering
$\alpha=\frac{1}{2}$
and
$\beta=1$
at the interval 
$0\leq x\leq 1$
in hyperbolic telegraph Eq. (\ref{qlm1}) with the following initial conditions:
\begin{equation*}
\hspace*{\fill}
\left\{
\begin{array}{ll}
u(x,0)=x^2 & 0\leq x\leq 1\\
u_t(x,0)=1, & 0\leq x\leq1
\end{array}
\right.
\hspace*{\fill}
\end{equation*}
we have
\begin{equation*}
\hspace*{\fill} \frac{\partial ^{2}u}{\partial t^{2}}+ \frac{\partial u}{\partial t} + u=\frac{\partial ^{2}u}{\partial x^{2}}+x^2+t-1,\hspace*{\fill}
\end{equation*}
Furthermore, the exact solution is
\begin{equation*}
\hspace*{\fill}u(x,t)=x^2+t.\hspace*{\fill}
\end{equation*}

Using the proposed method with layers and ten Legendre's polynomials as activation functions in the first layer and using initial points $t_1=0.6$ and $t_2=0.8$, the approximate solutions for $t_3 = 1$ are obtained. This example's numerical results on the interval $[0, 1]$ for $t_3=1$ have been depicted in figure \ref{fig4} with 200 training points together with the residual
function. Moreover, the relative $L_2$ errors in JDNN and the simple DNN for random training and testing data have been computed and presented in table \ref{1.2} with 2000 epochs.
Figures \ref{123} and \ref{1234} demonstrate the convergence of the loss function throughout the learning prosses. 
Table \ref{2.2} compares the results with the collocation points method presented by Dehghan et al. \cite{https://doi.org/10.1002/num.20306}.
\begin{figure}[H]
   \centering
	\begin{subfigure}[b]{0.4\textwidth}
  	\includegraphics[width=1\linewidth]{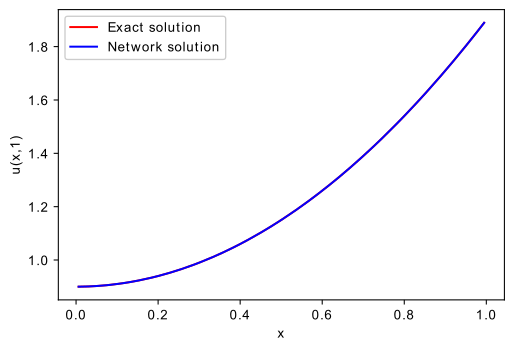}
      \caption{}
    \end{subfigure}
    \begin{subfigure}[b]{0.4\textwidth}
  	\includegraphics[width=1\linewidth]{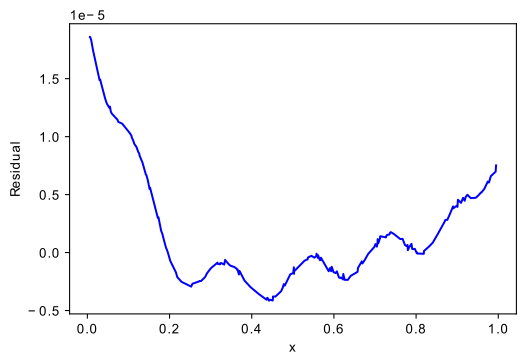}
      \caption{}
    \end{subfigure}
	\caption{(a) Exact and JDNN solutions for the example \ref{ex2} with training points in [0, 1] and (b) Absolute of the residual function.}
	 \label{fig4}
         
\end{figure}
\begin{figure}[H]
	\centering
	\begin{subfigure}[b]{0.4\textwidth}
  	\includegraphics[width=1.0\linewidth]{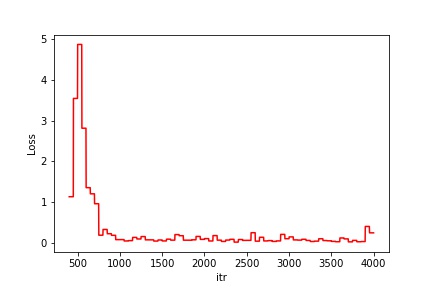}
      \caption{}
      \label{123}
    \end{subfigure}
    \begin{subfigure}[b]{0.4\textwidth}
  	\includegraphics[width=1.0\linewidth]{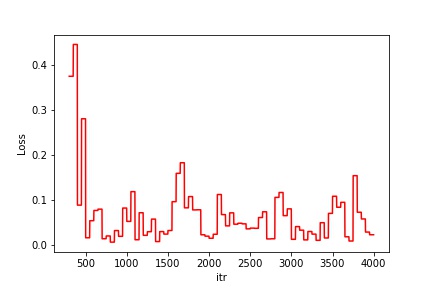}
      \caption{}
      \label{1234}
    \end{subfigure}
	\caption{Comparing the trend of Loss function changes in (a) the JDNN and (b) the simple DNN.}
\end{figure}
\begin{table*}[!hbtp]
	\centering
	\caption{Comparison of relative $L_2$ errors obtained in the JDNN and the simple DNN.}
	\begin{tabular}{|c|c|c|c|}
		\hline  \multirow{1}{*}{Method}& \multicolumn{1}{|c|}{Layers}& \multicolumn{1}{|c|}{Train data error}&  \multicolumn{1}{|c|}{Test data error}
		\\
		\hline JDNN  & $[1, 10, 20, 60, 80, 60, 20, 10, 1]$& $1.8\times 10^{-6}$ &  $2\times 10^{-6}$\\
		Simple DNN   & $[1, 10, 20, 60, 80, 60, 20, 10, 1]$& $4\times 10^{-6}$ &  $8\times 10^{-6}$ \\
		\hline
	\end{tabular}
	\label{1.2}
\end{table*}
\begin{table*}[!hbtp]
	\centering
	\caption{Comparison of the results obtained in the JDNN and the Collocation points method.}
	\begin{tabular}{|c|c|c|c|}
		\hline  \multirow{1}{*}{Method} &  \multicolumn{1}{|c|}{$L_{\infty}$-error} &  \multicolumn{1}{|c|}{$L_{2}$-error} &  \multicolumn{1}{|c|}{RMS}
		\\
		\hline JDNN  &  $2.0742\times 10^{-5}$&  $9.7068\times 10^{-5}$&  $6.8637\times 10^{-6}$\\
		Collocation points \cite{https://doi.org/10.1002/num.20306} &  $8.5573\times 10^{-5}$&  $6.1544\times 10^{-4}$&  $6.1239\times 10^{-5}$\\
		\hline
	\end{tabular}
	\label{2.2}
\end{table*}
\end{example}
\begin{example}\label{ex3}
Consider the second-order hyperbolic Eq. (\ref{qlm1n}) with $\alpha=1$ and $\beta=1$. The initial conditions are given by
\begin{equation*}
\hspace*{\fill}
\left\{
\begin{array}{ll}
u(x,y,0)=x_1^2+x_2^2 & 0\leq x_1, x_2\leq 1\\
u_t(x,y,0)=x_1^2+x_2^2+1, & 0\leq x_1, x_2 \leq1
\end{array}
\right.
\hspace*{\fill}
\end{equation*}
and the exact solution is
\begin{equation*}
\hspace*{\fill}u(x_1,x_2,t)=x_1^2+x_2^2+t,\quad   0\leq x_1, x_2 \leq1, \quad t>0.\hspace*{\fill}
\end{equation*}
The right hand side function is 
$f(x_1,x_2,t)=-2+x_1^2+x_2^2+t$
. We extract the boundary conditions from the exact solution.

By applying the 6-layered JDNN with the first 8 Chebyshev polynomials as activation functions for orbitary time $t_1=0.6$ and $t_2=0.8$, the approximate solutions for $t_3 = 1$ are calculated. The results and the residual function with 100 training points for solving the problem on the interval $[0, 1]\times[0, 1]$ are also shown in figure \ref{fig6}. Figure \ref{fig7} shows how the MSE loss function, defined in equation (\ref{loss}), converges during the training phase. In table \ref{1.3}, the results of the proposed networks are compared with the simple DNN. In addition, the results are compared to a meshless method proposed by Mehdi Dehghan and Ali Shokri in table \ref{2.3} \cite{https://doi.org/10.1002/num.20357}.

\begin{figure}[H]
	\centering
	\begin{subfigure}[b]{0.4\textwidth}
  	\includegraphics[width=1.05\linewidth]{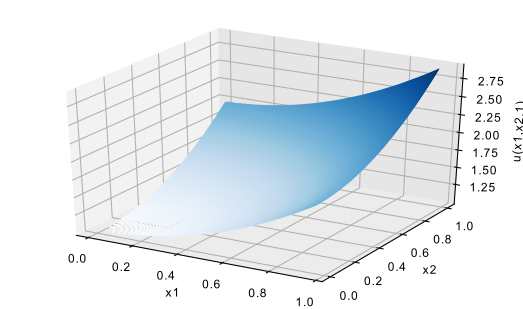}
      \caption{}
    \end{subfigure}
    \begin{subfigure}[b]{0.4\textwidth}
  	\includegraphics[width=1.05\linewidth]{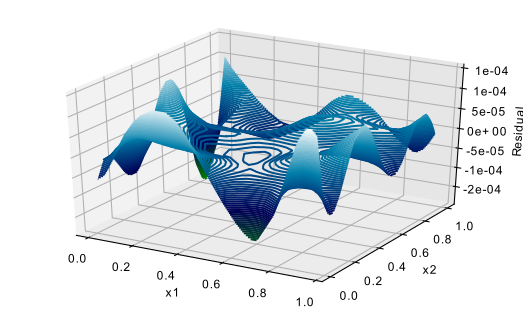}
      \caption{}
    \end{subfigure}
	\caption{(a) Exact and JDNN solutions for the example \ref{ex3} with training points and (b) Absolute of the residual function.}
         \label{fig6}
\end{figure}
\begin{figure}[H]
	\centering
	\begin{subfigure}[b]{0.4\textwidth}
  	\includegraphics[width=1.0\linewidth]{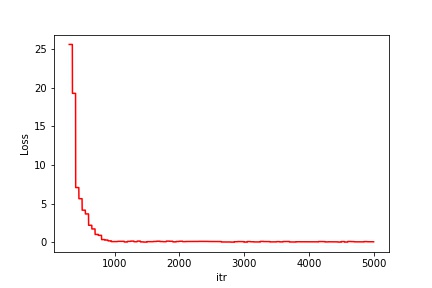}
      \caption{}
    \end{subfigure}
    \begin{subfigure}[b]{0.4\textwidth}
  	\includegraphics[width=1.0\linewidth]{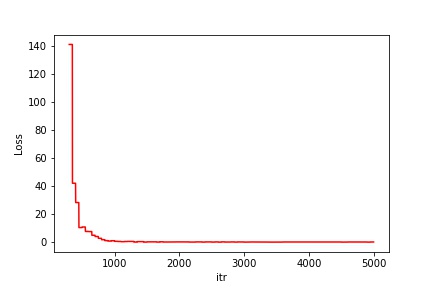}
      \caption{}
    \end{subfigure}
	\caption{Comparing the trend of Loss function changes in (a) the JDNN and (b) the simple DNN.}
	\label{fig7}
\end{figure}
\begin{table*}[!hbtp]
	\centering
	\caption{Comparison of relative $L_2$ errors obtained in the JDNN and simple DNN.}
	\begin{tabular}{|c|c|c|c|}
		\hline  \multirow{1}{*}{Method}& \multicolumn{1}{|c|}{Layers}& \multicolumn{1}{|c|}{Train data error}&  \multicolumn{1}{|c|}{Test data error}
		\\
		\hline JDNN  & $[2, 8, 16, 32, 64, 32, 16, 8, 1]$& $8.7\times 10^{-6}$ &  $1.6\times 10^{-5}$\\
		Simple DNN  & $[2, 8, 16, 32, 64, 32, 16, 8, 1]$& $2.4\times 10^{-5}$ &  $2.9\times 10^{-5}$ \\
		\hline
	\end{tabular}
	\label{1.3}
\end{table*}
\begin{table*}[!hbtp]
	\centering
	\caption{Comparison of the results obtained in the JDNN and Meshless Method.}
	\begin{tabular}{|c|c|c|c|}
		\hline  \multirow{1}{*}{Method} &  \multicolumn{1}{|c|}{$L_{\infty}$-error} &  \multicolumn{1}{|c|}{$L_{2}$-error} &  \multicolumn{1}{|c|}{RMS}
		\\
		\hline JDNN  &  $8.5818\times 10^{-5}$&  $9.1080\times 10^{-4}$&  $2.3720\times 10^{-5}$\\
		Meshless Method \cite{https://doi.org/10.1002/num.20357}  &  $1.8056\times 10^{-4}$&  $1.2439\times 10^{-3}$&  $1.1309\times 10^{-4}$\\
		\hline
	\end{tabular}
	\label{2.3}
\end{table*}
\end{example}
\begin{example}\label{ex4}
Consider the second-order hyperbolic Eq. (\ref{qlm1n}) with $\alpha=1$ and $\beta=1$. The initial conditions are given by
\begin{equation*}
\hspace*{\fill}
\left\{
\begin{array}{ll}
u(x_1,x_2,0)=sin(x_1) sin(x_2) & 0\leq x_1, x_2\leq 1\\
u_t(x_1,x_2,0)=0, & 0\leq x_1, x_2 \leq1
\end{array}
\right.
\hspace*{\fill}
\end{equation*}
and the exact solution is
\begin{equation*}
\hspace*{\fill}u(x_1,x_2,t)=cos(t) sin(x_1) sin(x_2),\quad   0\leq x_1, x_2 \leq1,\quad t>0.\hspace*{\fill}
\end{equation*}
The right hand side function is 
$f(x_1,x_2,t)=2sin(x_1) sin(x_2) (cos(t)-sin(t))$
. We extract the boundary conditions from the exact solution.

We Use the 6-layered JDNN with the first 4 Chebyshev polynomials as activation functions; this example's numerical results on the interval $[0, 1]\times[0, 1]$ for $t_3=1$ have been obtained, depicted in figure \ref{fig8} with 100 training points together with the residual
function. In table \ref{1.4}, the results of the proposed network are compared with the simple DNN. Also, in table \ref{2.4}, the results are compared with the Generalized Finite Difference Method presented by F. Ureña et al., reflecting the proposed method's good performance \cite{ar44}.
Figure \ref{fig9} compares the loss function values for the JDNN and the Simple DNN models for the numerical simulation of problem 4. 
\begin{figure}[H]
	\centering
	\begin{subfigure}[b]{0.4\textwidth}
  	\includegraphics[width=1\linewidth]{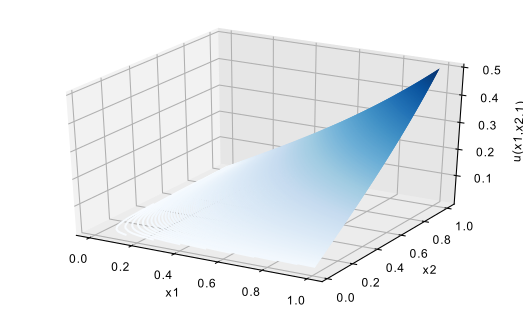}
      \caption{}
    \end{subfigure}
    \begin{subfigure}[b]{0.4\textwidth}
  	\includegraphics[width=1.06\linewidth]{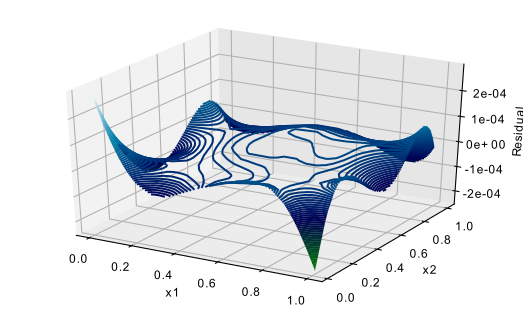}
      \caption{}
    \end{subfigure}
	\caption{(a) Exact and JDNN solutions for the example \ref{ex4} with training points and (b) Absolute of the residual function.}
          \label{fig8}
\end{figure}
\begin{figure}[H]
	\centering
	\begin{subfigure}[b]{0.4\textwidth}
  	\includegraphics[width=1.1\linewidth]{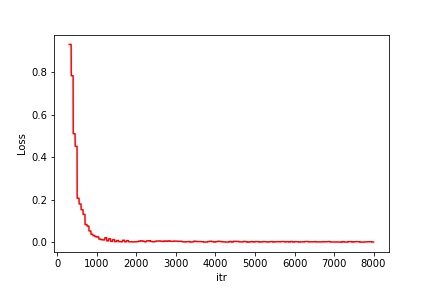}
      \caption{}
    \end{subfigure}
    \begin{subfigure}[b]{0.4\textwidth}
  	\includegraphics[width=1.1\linewidth]{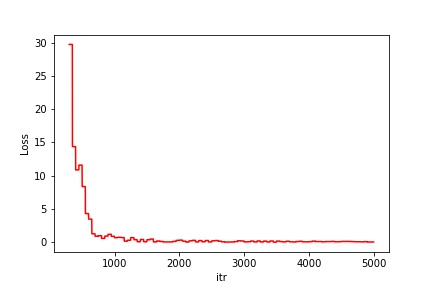}
      \caption{}
    \end{subfigure}
	\caption{Comparing the trend of Loss function changes in (a) the JDNN and (b) the simple DNN.}
	\label{fig9}
\end{figure}
\begin{table*}[!hbtp]
	\centering
	\caption{Comparison of relative $L_2$ errors obtained in the JDNN and the simple DNN.}
	\begin{tabular}{|c|c|c|c|}
		\hline  \multirow{1}{*}{Method}& \multicolumn{1}{|c|}{Layers}&  \multicolumn{1}{|c|}{Train data error}&  \multicolumn{1}{|c|}{Test data error}
		\\
		\hline JDNN & $[2, 4, 8, 16, 32, 32, 16, 8, 1]$ & $5.4\times 10^{-4}$ &  $6.3\times 10^{-4}$\\
		Simple DNN & $[2, 4, 8, 16, 32, 32, 16, 8, 1]$ & $7.58\times 10^{-4}$ &  $1.1\times 10^{-3}$ \\
		\hline
	\end{tabular}
	\label{1.4}
\end{table*}
\begin{table*}[!hbtp]
	\centering
	\caption{Comparison of the results obtained in the JDNN and GFDM.}
	\begin{tabular}{|c|c|c|}
		\hline  \multirow{1}{*}{Method} &  \multicolumn{1}{|c|}{$L_{\infty}$-error} &  \multicolumn{1}{|c|}{$L_{2}$-error} 
		\\
		\hline JDNN  &  $3.10\times 10^{-4}$&  $2.3\times 10^{-4}$\\
		GFDM \cite{ar44}  &  $5.18\times 10^{-4}$&  $2.44\times 10^{-4}$\\
		\hline
	\end{tabular}
	\label{2.4}
\end{table*}
\end{example}
\section{Conclusion}
In this article, we Have developed Jacobi deep neural network (JDNN) for solving 1-D and 2-D telegraph equations. The purpose of introducing this network is to take advantage of the collocation method and deep learning. We Have used the finite difference method to discretize the residual form to utilize the proposed method to solve the telegraph equation. JDNN consists of two parts, a feedforward neural network for approximating the solution of the equation and a second part for calculating the residual form of the
equation, which is discretized using the finite difference method. We have solved the telegraph equation once with the JDNN and once with a similar neural network that does not use the orthogonal polynomials activation functions. All of its activation functions are Tanh (Simple DNN). A comparison of this network's results shows the proposed method's stability. It can also be seen that the JDNN has better results than a simple DNN. In the architecture of this network, we used two kinds of Jacobi polynomials, Chebyshev and Legendre polynomials, as activation functions to improve the results. The numerical results
applying the training and testing points show the high accuracy of our proposed model. In addition, the JDNN can be used to solve other PDEs of higher degrees, similar to the telegraph equation.


\bibliography{main}

\begin{thebibliography}{10}
\expandafter\ifx\csname url\endcsname\relax
  \def\url#1{\texttt{#1}}\fi
\expandafter\ifx\csname urlprefix\endcsname\relax\def\urlprefix{URL }\fi
\expandafter\ifx\csname href\endcsname\relax
  \def\href#1#2{#2} \def\path#1{#1}\fi

\bibitem{par}
A.~Kolesov, N.~K. Rozov, {Parametric} excitation of high-mode oscillations for
  a non-linear telegraph equation, Sbornik: Mathematics 191(8) (2000)
  1147–69.

\bibitem{ind}
A.~C. Metaxas, R.~J. Meredith, Industrial {Microwave} {Heating}, no.~4, IET,
  1983.

\bibitem{ar44}
F.~Prieto, L.~Gavete, J.~Benito, A.~García, A.~Vargas, Solving the telegraph
  equation in {2-D} and {3-D} using {Generalized} {Finite} {Difference}
  {Method} {(GFDM)}, Engineering Analysis with Boundary Elements 112 (2019)
  13--24.

\bibitem{articl3e}
Y.~Youssri, W.~Abd-Elhameed, Numerical {Spectral} {Legendre-Galerkin}
  {Algorithm} for {Solving} {Time} {Fractional} {Telegraph} {Equation}, Rom. J.
  Phys 63~(107) (2018) 1--16.

\bibitem{irk2019numerical}
D.~Irk, E.~Kirli, Numerical solution of the homogeneous telegraph equation by
  using galerkin finite element method, International Conference on
  Computational Mathematics and Engineering Sciences (2019) 209--217.

\bibitem{doi:10.2514/3.559}
J.~Chai, H.~Lee, S.V.Patankar, Finite volume method for radiation heat
  transfer, Journal of Thermophysics and Heat Transfer 8~(3) (1994) 419--425.

\bibitem{Wrobel2003BoundaryEM}
L.~Wrobel, A.~Kassab, {{Boundary Element Method}, {Volume} 1: Applications in
  {Thermo-Fluids} and Acoustics}, Applied Mechanics Reviews 56~(2) (2003)
  B17--B17.

\bibitem{https://doi.org/10}
T.~Belytschko, Y.~Y. Lu, L.~Gu, Element-free {Galerkin} methods, International
  Journal for Numerical Methods in Engineering 37~(2) (1994) 229--256.

\bibitem{article3}
S.~Atluri, T.~Zhu, The {Meshless} {Local} {Petrov-Galerkin (MLPG)} approach for
  solving problems in elasto-statics, Computational Mechanics 25 (2000)
  169--179.

\bibitem{article1}
A.~Kumar, A.~Bhardwaj, S.~Dubey, A local meshless method to approximate the
  time-fractional telegraph equation, Engineering with Computers (2020) 1--16.

\bibitem{artJ}
R.~Jiwari, {Lagrange} interpolation and modified cubic {B-spline} differential
  quadrature methods for solving hyperbolic partial differential equations with
  {Dirichlet} and {Neumann} boundary conditions, Computer Physics
  Communications 193 (2015) 55--65.

\bibitem{arO}
{\"O}.~Oru{\c{c}}, A numerical procedure based on {Hermite} wavelets for
  two-dimensional hyperbolic telegraph equation, Engineering with Computers
  34~(4) (2018) 741--755.

\bibitem{artiji}
R.~Jiwari, S.~Pandit, R.~Mittal, A {Differential} {Quadrature} {Algorithm} for
  the {Numerical} {Solution} of the {Second-Order} {One} {Dimensional}
  {Hyperbolic} {Telegraph} {Equation}, International Journal of Non-linear
  Sciences 13 (2012) 259--266.

\bibitem{hajimohammadi2021legendre}
Z.~Hajimohammadi, K.~Parand, A.~Ghodsi, {Legendre} {Deep} {Neural} {Network}
  {(LDNN)} and its application for approximation of nonlinear {Volterra}
  {Fredholm} {Hammerstein} integral equations, arXiv preprint arXiv:2106.14320.

\bibitem{6224185}
S.~Mehrkanoon, T.~Falck, J.~A.~K. Suykens, Approximate {Solutions} to {Ordinary
  Differential Equations} {Using Least Squares Support Vector Machines}, IEEE
  Transactions on Neural Networks and Learning Systems 23~(9) (2012)
  1356--1367.

\bibitem{GoodBengCour16}
I.~Goodfellow, Y.~Bengio, A.~Courville, Deep Learning, MIT press, 2016.

\bibitem{MEHR}
S.~Mehrkanoon, J.~Suykens, Learning solutions to partial differential equations
  using {LS-SVM}, Neurocomputing 159 (2015) 105--116.

\bibitem{SIRI}
J.~Sirignano, K.~Spiliopoulos, {DGM}: A deep learning algorithm for solving
  partial differential equations, Journal of Computational Physics 375.

\bibitem{RAISI}
M.~Raissi, P.~Perdikaris, G.~Karniadakis, {Physics-Informed} {Neural}
  {Networks}: A deep learning framework for solving forward and inverse
  problems involving nonlinear partial differential equations, Journal of
  Computational Physics 378 (2019) 686--707.

\bibitem{padierna2018novel}
L.~C. Padierna, M.~Carpio, A.~Rojas-Dom{\'\i}nguez, H.~Puga, H.~Fraire, A novel
  formulation of orthogonal polynomial kernel functions for svm classifiers:
  the gegenbauer family, Pattern Recognition 84 (2018) 211--225.

\bibitem{wavelet2}
L.~Zhang, W.~Zhou, L.~Jiao, Wavelet support vector machine, IEEE Transactions
  on Systems, Man, and Cybernetics, Part B (Cybernetics) 34~(1) (2004) 34--39.

\bibitem{wavelet}
M.~R. Gupta, N.~P. Jacobson, Wavelet principal component analysis and its
  application to hyperspectral images, in: 2006 International Conference on
  Image Processing, 2006, pp. 1585--1588.

\bibitem{HAJI}
Z.~Hajimohammadi, K.~Parand, Numerical learning approximation of
  time-fractional sub diffusion model on a semi-infinite domain, Chaos,
  Solitons and Fractals 142 (2021) 110435.

\bibitem{parand2021parallel}
K.~Parand, A.~A. Aghaei, M.~Jani, A.~Ghodsi, Parallel ls-svm for the numerical
  simulation of fractional volterra’s population model, Alexandria
  Engineering Journal 60~(6) (2021) 5637--5647.

\bibitem{parand2}
A.~G. Khoee, K.~M. Mohammadi, M.~Jani, K.~Parand, {A least squares support
  vector regression for anisotropic diffusion filtering}, arXiv preprint
  arXiv:2202.00595.

\bibitem{HADI}
A.~Hadian-Rasanan, D.~Rahmati, S.~Gorgin, K.~Parand, A {Single} {Layer}
  {Fractional} {Orthogonal} {Neural} {Network} for {Solving} {Various} {Types}
  of {Lane-Emden Equation}, New Astronomy 75 (2020) 101307.

\bibitem{hajim}
Z.~Hajimohammadi, F.~Baharifard, {Fractional Chebyshev deep neural network
  (FCDNN) for solving differential models}, Chaos, Solitons \& Fractals 153
  (2021) 111530.

\bibitem{chak}
S.~Chakraverty, S.~Mall, {Artificial} {Neural} {Networks} for {Engineers} and
  {Scientists}: {Solving} {Ordinary} {Differential} {Equations}, CRC Press,
  2017.

\bibitem{doi:10.1063/1.369258}
P.~M. Jordan, A.~Puri, {Digital} signal propagation in dispersive media,
  Journal of Applied Physics 85~(3) (1999) 1273--1282.

\bibitem{article22}
J.~Banasiak, J.~Mika, Singular perturbed telegraph equations with applications
  in random walk theory, Journal of Applied Mathematics and Stochastic Analysis
  11~(1) (1998) 9--28.

\bibitem{Roussy1998FoundationsAI}
G.~Roussy, J.~Pearce, {Foundations And Industrial Applications Of Microwaves
  And Radio Frequency Fields. Physical And Chemical Processes}, Proceedings of
  the 6th International Conference on Optimization of Electrical and Electronic
  Equipments 1 (1998) 115--116.

\bibitem{article44}
D.~Rostamy, M.~Emamjomeh, S.~Abbasbandy, A meshless technique based on the
  pseudospectral radial basis functions method for solving the two-dimensional
  hyperbolic telegraph equation, The European Physical Journal Plus 132~(6)
  (2017) 1--11.

\bibitem{MA2016236}
W.~Ma, B.~Zhang, H.~Ma, A meshless collocation approach with barycentric
  rational interpolation for two-dimensional hyperbolic telegraph equation,
  Applied Mathematics and Computation 279 (2016) 236--248.

\bibitem{doi:10.1080/00207160211918}
R.~K. Mohanty, M.~K. Jain, A.~Urvashi, An {Unconditionally} {Stable} {ADI}
  {Method} for the linear {Hyperbolic} {Equation} in {Three} {Space}
  {Dimensions}, International Journal of Computer Mathematics 79~(1) (2002)
  133--142.

\bibitem{funaro2008polynomial}
D.~Funaro, {Polynomial Approximation of Differential Equations}, Vol.~8,
  Springer Science \& Business Media, 2008.

\bibitem{johnson2012numerical}
C.~Johnson, {Numerical} {Solution} of {Partial} {Differential} {Equations} by
  the {Finite} {Element} {Method}, Courier Corporation, 2012.

\bibitem{adam}
D.~P. Kingma, J.~Ba, Adam: a method for stochastic opoimization (2017).

\bibitem{bfgs}
C.~Zhu, R.~H. Byrd, P.~Lu, J.~Nocedal, {Algorithm 778: L-BFGS-B: Fortran
  subroutines for large-scale bound-constrained optimization}, ACM Transactions
  on mathematical software (TOMS) 23~(4) (1997) 550--560.

\bibitem{MOMANI20051126}
S.~Momani, Analytic and approximate solutions of the space- and time-fractional
  telegraph equations, Applied Mathematics and Computation 170~(2) (2005)
  1126--1134.

\bibitem{https://doi.org/10.1002/num.20306}
M.~Dehghan, A.~Shokri, {A Numerical Method for Solving the Hyperbolic Telegraph
  Equation}, Numerical Methods for Partial Differential Equations 24~(4) (2008)
  1080--1093.

\bibitem{https://doi.org/10.1002/num.20357}
M.~Dehghan, A.~Shokri, A {Meshless} method for numerical solution of a linear
  hyperbolic equation with variable coefficients in two space dimensions,
  Numerical Methods for Partial Differential Equations 25~(2) (2009) 494--506.

\end{thebibliography}

\end{document}